\documentclass[10pt]{article}
\usepackage{graphicx}
\usepackage{amsmath}
\usepackage{amsfonts}
\usepackage{amssymb}

\def\d{\mathrm d}
\def\e{\mathrm e}

\def\limfdd{\renewcommand{\arraystretch}{0.5}
\begin{array}[t]{c}
\stackrel{\rm fdd}{\longrightarrow} \\
\end{array}\renewcommand{\arraystretch}{1}}

\def\i{\mathrm i}

\newtheorem{thm}{Theorem}[section]

\newtheorem{prop}[thm]{Proposition}

\newtheorem{defn}[thm]{Definition}

\newtheorem{rem}[thm]{Remark}

\numberwithin{equation}{section}

\newcommand{\R}{\mathbb R}

\newcommand{\E}{\mathbb E}

\newcommand{\Z}{\mathbb Z}

\def\eqfdd{\renewcommand{\arraystretch}{0.5}
\begin{array}[t]{c}
\stackrel{\rm fdd}{=} \\
\end{array}\renewcommand{\arraystretch}{1}}

\def\eqd{\renewcommand{\arraystretch}{0.5}
\begin{array}[t]{c}
\stackrel{\rm d}{=} \\
\end{array}\renewcommand{\arraystretch}{1}}

\newcommand{\nn}{\nonumber}
\newcommand{\noi}{\noindent}

\def\limfdd{\renewcommand{\arraystretch}{0.5}
\begin{array}[t]{c}
\stackrel{\rm fdd}{\longrightarrow} \\
\end{array}\renewcommand{\arraystretch}{1}}

\def\limfdd{\renewcommand{\arraystretch}{0.5}
\begin{array}[t]{c}
\stackrel{\rm fdd}{\longrightarrow} \\
\end{array}\renewcommand{\arraystretch}{1}}

\begin{document}

\title{Tempered fractional Brownian and stable motions \\of  second kind }

\author{Farzad Sabzikar$^{1}$  \  and \  Donatas Surgailis$^2$}
\date{\today \\  \small
\vskip.2cm
$^1$Iowa State University \ and \
$^2$Vilnius University}

\maketitle

\begin{abstract}  Meerschaert and Sabzikar \cite{Meerschaertsabzikar}, \cite{Meerschaertsabzikar2}
introduced tempered fractional Brownian/stable motion (TFBM/TFSM) 
by including an exponential  tempering factor
in the moving  average representation of FBM/FSM.  
The present paper discusses another tempered
version of FBM/FSM, termed tempered fractional Brownian/stable motion of second kind (TFBM II/TFSM II).
We prove that TFBM/TFSM and TFBM II/TFSM II are different processes. Particularly, large time properties
of TFBM II/TFSM II are similar to those of FBM/FSM and are in deep contrast to large time properties
of  TFBM/TFSM.
\end{abstract}

\smallskip
{\small

\noi {\it Keywords:} 
tempered fractional Brownian/stable motion; tempered fractional Brownian/stable noise;
tempered fractional integration;
local and global self-similarity

}

\section{Introduction}


Meerschaert and Sabzikar \cite{Meerschaertsabzikar2}
introduced tempered fractional stable motion (TFSM) $ Z_{H,\alpha,\lambda} =  \{Z_{H,\alpha,\lambda}(t), t \in \R \} $
for $0< \alpha \le 2, H> 0, \lambda >0$
as stochastic integral
\begin{eqnarray}\label{TFSM}
Z_{H,\alpha,\lambda}(t)
&:=&\int_{\R}\Big( (t-y)_+^{H - \frac{1}{\alpha}} \e^{-\lambda (t-y)_+} - (-y)_+^{H - \frac{1}{\alpha}} \e^{-\lambda (-y)_+}
\Big) M_\alpha (\d y),
\end{eqnarray}
with respect to $\alpha$-stable  L\'evy process $M_\alpha $.
A particular case of TFSM termed the {\it tempered fractional Brownian motion} (TFBM)
corresponding to $\alpha =2$ and $M_2 = B$ (a standard Brownian motion) was studied
in Meerschaert and Sabzikar \cite{Meerschaertsabzikar}.
Note that
for $\lambda  = 0$ (and $H \in (0,1)$) TFSM/TFBM agree with
fractional stable/Brownian motion (FSM/FBM), see \cite{SamorodnitskyTaqqu}. The role of the tempering by exponential factor
in \eqref{TFSM}
manifests in the dependence properties of the increment process
$Y_{H, \alpha, \lambda} =
\{Y_{H, \alpha, \lambda}(t) := Z_{H, \alpha, \lambda}(t+1) - Z_{H, \alpha, \lambda}(t), t \in \Z \} $
called {\it tempered fractional stable noise} (TFSN) and {\it tempered fractional Gaussian noise} (TFGN) in the Gaussian
case $\alpha = 2$. In particular, for small $\lambda >0$
the autocovariance function of TFGN
closely resembles that of fractional Gaussian noise
(FGN) on an intermediate scale, but then it eventually falls off exponentially. On the other hand,
the spectral density of TFGN vanishes at the origin for all $H >0$ exhibiting a strong
anti-persistent behavior, see \cite{Meerschaertsabzikar}.

In this paper we study a closely related but different tempered process
called  {\it tempered fractional Brownian/stable motion of second kind} (TFBM II/TFSM II)
which is defined by replacing the integrand in \eqref{TFSM} by
\begin{eqnarray}\label{hdef0}
h_{H,\alpha,\lambda}(t;y)&:=&(t-y)_+^{H - \frac{1}{\alpha}} \e^{-\lambda (t-y)_+} - (-y)_+^{H - \frac{1}{\alpha}} \e^{-\lambda (-y)_+}
+ \lambda \int_{0}^{t} (s-y)_{+}^{H-\frac{1}{\alpha}} \e^{-\lambda(s-y)_{+}}\ \d s, \qquad y \in \R. \nn
\end{eqnarray}
The corresponding $\alpha$-stable process, denoted by $Z^{I\!I}_{H,\alpha,\lambda} =
\{ Z^{I\!I}_{H,\alpha,\lambda}(t), t \in \R \}$ is defined for all $H>0, 1< \alpha \le 2, \lambda >0$ and has
stationary increments similarly as $Z_{H,\alpha,\lambda}$. The change of the integrand results
in a drastic change of large-time behavior of the increment process
\begin{equation} \label{eq:TFSN definition}
Y^{I\!I}_{H, \alpha, \lambda} =
\{Y^{I\!I}_{H, \alpha, \lambda}(t) := Z^{I\!I}_{H, \alpha, \lambda}(t+1) - Z^{I\!I}_{H, \alpha, \lambda}(t), t \in \Z \}
\end{equation}
called {\it tempered fractional Gaussian/stable noise of second kind} (TFGN II/TFSN II).
The spectral density of TFGN II
$ Y^{I\!I}_{H,2, \lambda}$ decays as a power function for frequencies $|\omega| > \lambda $  but
remains bounded and separated from zero near zero frequency,
making
TFGN II a realistic model in turbulence and other applied areas. See Figures 1, 2 and Remark \ref{kolm}.

One of the main motivation for our introducing and studying  is the fact that these processes appear
as the limits of the partial sums process
of tempered stationary linear processes
with discrete time
and small tempering parameter $\lambda_N \sim \lambda/N $ tending to zero together with the sample size $N$.
This problem is discussed in a  parallel paper Sabzikar and Surgailis \cite{SabzikarSurgailis} where we prove that
the limit behavior of such partial process essentially  depends on how fast the tempering parameter
tends to zero, resulting in different limits in the strongly tempered ($\lim_{N \to \infty} \lambda_N/N = 0$),
weakly tempered ($\lim_{N \to \infty} \lambda_N/N = \infty$), and moderately tempered ($\lim_{N \to \infty} \lambda_N/N \in (0,\infty)$)
situations.

Let us describe the main results of this paper.
Section 2 provides the basic definitions and properties of
TFBM II/TFSM II. The latter include the spectral representation and
the covariance function of TFBM II, relation to tempered fractional calculus
(see \cite{MeerschaertsabzikarSPA}), and the relation between TFSM and TFSM II.
Theorem \ref{lem:scaling}
establishes local and global asymptotic self-similarity of TFSM and TFSM II.
It shows that TFSM and TFSM II are very different processes; indeed, the former process
is stochastically bounded and the latter is stochastically unbounded as $t \to \infty$.
Section 3 discusses the dependence properties of stationary processes
TFSN II and TFGN II. We obtain the asymptotic behavior of the bivariate characteristic function of
TFSN II which can be compared to the corresponding results for TFSN
in  Meerschaert and Sabzikar in \cite{Meerschaertsabzikar2} and  for FSN
in Astrauskas et al.  \cite{Astrauskas}.

In what follows, $C$ denotes generic constants 
which may be different at different locations.
We write $
\limfdd $ and  $   \eqfdd $  
 for weak convergence and equality
of finite-dimensional
distributions, respectively. $\R_+ := (0, \infty),
(x)_\pm := \max (\pm x, 0),  x \in \R,  \, \int := \int_\R$.

\section{Definition and properties of TFSM II and TFBM II
}\label{sec2}

For $0< \alpha \le 2$,
let  $\{M_\alpha(t)\}_{t\in \R}$ be an $\alpha$-stable 
L\'evy process with stationary independent increments and characteristic
function
\begin{equation}\label{Mcf}
\E \e^{\i \theta M_\alpha (t)}  = \e^{ - \sigma^\alpha |\theta|^\alpha t
\big(1 - \i \beta \tan (\pi \alpha/2) {\rm sign}(\theta)\big)},
\qquad \theta \in \R,
\end{equation}
where $\sigma >0 $ and $ \beta \in [-1,1]$
are  the scale and skewness parameters, respectively. For $\alpha =2$,
$M_2(t) = \sqrt{2} \sigma B(t)$, where $B $ is a standard Brownian motion with variance $\E B^2(t) = t$. Stochastic
integral $I_\alpha(f) \equiv \int f(x) M_\alpha (\d x)$ is defined for any $ f \in L^\alpha (\R)$
as $\alpha$-stable
random variable with
characteristic function
\begin{equation} \label{Ialpha}
\E \e^{ \i \theta I_\alpha (f)}  = \exp \{  - \sigma^\alpha |\theta|^\alpha
\int |f(x)|^\alpha
\big(1 - \i \beta \tan (\pi \alpha/2) {\rm sign}(\theta f(x))\big)  \d x \},
\quad \theta \in \R.
\end{equation}
see e.g. \cite[Chapter 3]{SamorodnitskyTaqqu}.

Note the function  $y \mapsto h_{H,\alpha,\lambda}(t;y): \R \to \R $ in \eqref{hdef0}  satisfies  $h_{H, \alpha,\lambda}(t;\cdot) \in L^\alpha (\R)  $ for any $t \in \R $ and any  $\lambda >0,   1< \alpha \le 2,  H > 0$  and also for  $\lambda = 0, 0< \alpha \le 2,
H \in (0,1)$.
We will use the following integral representation of \eqref{hdef0}: \\
For $H > \frac{1}{\alpha}$:
\begin{eqnarray}\label{hdef1}
&h_{H,\alpha,\lambda}(t;y) = (H - \frac{1}{\alpha})
\int_0^t (s-y)_+^{H - \frac{1}{\alpha} -1} \e^{-\lambda (s-y)_+} \d s.
\end{eqnarray}
For $0< H < \frac{1}{\alpha}$:
\begin{eqnarray} \label{hdef2}
&h_{H,\alpha,\lambda}(t;y)=(H - \frac{1}{\alpha})
\begin{cases}\int_0^t (s-y)_+^{H - \frac{1}{\alpha} -1} \e^{-\lambda (s-y)_+} \d s, &y < 0, \\
-\int_t^\infty (s-y)_+^{H - \frac{1}{\alpha} -1} \e^{-\lambda (s-y)_+} \d s + \lambda^{\frac{1}{\alpha} -H}
\Gamma(H - \frac{1}{\alpha})
, &y \ge 0.
\end{cases}
\end{eqnarray}

\begin{defn}\label{defTHP}
{\rm
{Let $M_\alpha $ be $\alpha$-stable L\'evy process in \eqref{Mcf}, $1< \alpha \le 2 $ and
$H >0, \, \lambda > 0$. The stochastic process
\begin{eqnarray}\label{TFSM2}
Z^{I\!I}_{H,\alpha,\lambda}(t)&:=&\int  h_{H,\alpha,\lambda}(t;y)
M_\alpha (\d y), \qquad t \in \R
\end{eqnarray}
will be called
{\it tempered fractional stable motion of  second kind} (TFSM II).
A particular case of \eqref{TFSM2} corresponding to $\alpha = 2$
\begin{eqnarray}\label{TFBM2}
B^{I\!I}_{H,\lambda}(t)
&:=&\frac{1}{\Gamma(H+\frac{1}{2})}\int  h_{H,2,\lambda}(t;y) B(\d y)
\end{eqnarray}
will be called  {\it tempered  fractional Brownian motion of second kind} (TFBM II). } }
\end{defn}

\begin{rem} {\rm TFSM II can be also defined for $0 < \alpha \le 1 $, however
this definition uses a different kernel from \eqref{hdef0}
for $ 0 < H < \frac{1}{\alpha} -1.  $  The reason is that the tempered fractional derivative ${\mathbb{D}}^{\kappa,\lambda}_{\pm}$ in
\eqref{D} takes a different form for $\kappa > 1 $, see
(\cite{MeerschaertsabzikarSPA}, Remark 2.12).
}
\end{rem}

\begin{rem} \label{remstable} {\rm
For $\lambda >0$ and the same L\'evy process $M_\alpha $,
TFSM in \eqref{TFSM}
and TFSM II in \eqref{TFBM2}
are related as
\begin{equation} \label{ZX}
Z^{I\!I}_{H,\alpha,\lambda}(t) = Z_{H,\alpha,\lambda}(t) + \lambda C_{H,\alpha,\lambda}(t),
\end{equation}
where
\begin{equation}\label{CX}
C_{H,\alpha,\lambda}(t) := \int_0^t  Z_{H,\alpha,\lambda}(s) \d s + t C^0_{H,\alpha,\lambda}
\end{equation}
and
\begin{equation} \label{C0}
C^0_{H,\alpha,\lambda}  := \int  (-y)_+^{H - \frac{1}{\alpha}} \e^{-\lambda (-y)_+} M_\alpha (\d y)
\end{equation}
is a well-defined $\alpha$-stable r.v.  Relation \eqref{ZX} shows that $Z^{I\!I}_{H,\alpha,\lambda}$ and
$Z_{H,\alpha,\lambda}$ have similar path properties, particularly, path properties
of  $Z^{I\!I}_{H,\alpha,\lambda}$ can be derived from the path properties of $Z_{H,\alpha,\lambda}$ studied
in \cite{Meerschaertsabzikar2}.}
\end{rem}

Recall that
for any function in $L^{p}(\R)$, where $1\leq p<\infty$,
the (positive and negative) tempered fractional integrals are defined by
\begin{equation}\label{I}
{\mathbb I}^{\kappa,\lambda}_{\pm} f(y) :=\frac{1}{\Gamma(\kappa)}\int f(s)(y-s)_{\pm}^{\kappa-1}\e^{-\lambda(y-s)_{\pm}}\d s,
\quad \kappa >0
\end{equation}
and tempered fractional derivatives are
\begin{equation}\label{D}
{\mathbb{D}}^{\kappa,\lambda}_{\pm}f(y) := {\lambda}^{\kappa}f(y)+
\frac{\kappa}{\Gamma(1-\kappa)}\int \frac{f(y)-f(s)}{(y-s)_\pm^{\kappa+1}}
\,\e^{-\lambda(y-s)_\pm}\d s,  \quad 0< \kappa <1.
\end{equation}

The following proposition shows that TFSM II can be written as
a stochastic integral of tempered fractional integral (derivative) of the indicator function
of the interval $[0,t]$. In contrast, the corresponding expression for  TFSM in
\cite{MeerschaertsabzikarSPA} is more complicated and involves a linear combination
of tempered integrals and derivatives of the indicator function.


\begin{prop}\label{lem:connection with fractional calculus}
{\rm Let $H>0$, $1<\alpha\leq 2$, and $\lambda>0$. Then
\begin{eqnarray}
{\mathbb I}^{H - \frac{1}{\alpha},\lambda}_{-} {\bf 1}_{[0,t]}(y)
&=&\Gamma (H + 1 - \frac{1}{\alpha})^{-1} h_{H,  \alpha, \lambda}(t;y),
\quad  H>\frac{1}{\alpha}, \label{I1} \\
{\mathbb D}^{\frac{1}{\alpha}-H,\lambda}_{-} {\bf 1}_{[0,t]}(y)
&=&\Gamma (H + 1 - \frac{1}{\alpha})^{-1} h_{H,  \alpha, \lambda}(t;y),
\quad 0 < H< \frac{1}{\alpha}.   \label{D1}
\end{eqnarray}}
Thus, TFSM II for $1 < \alpha \le 2, H>0 $
can be defined as
\begin{equation}\label{eq:connection with fractional calculus}
Z^{I\!I}_{H,\alpha,\lambda}(t) =
 \Gamma (H + 1 - \frac{1}{\alpha}) \begin{cases}
\int   {\mathbb I}^{H - \frac{1}{\alpha},\lambda}_{-} {\bf 1}_{[0,t]}(y) M_\alpha (\d y), &  H>\frac{1}{\alpha}, \\
\int   {\mathbb D}^{\frac{1}{\alpha} - H,\lambda}_{-} {\bf 1}_{[0,t]}(y) M_\alpha (\d y),  &  0<  H <  \frac{1}{\alpha}.
\end{cases}
\end{equation}
\end{prop}

\noi {\it Proof.} Apply the tempered fractional operator ${\mathbb I}^{\kappa,\lambda}_{-}f$ in \eqref{I}
for $\kappa=H-\frac{1}{\alpha}$, $\lambda>0$ and
$f :={\bf 1}_{[0,t]}$ to see that
\begin{equation*}
\begin{split}
{\mathbb I}^{H-\frac{1}{\alpha},\lambda}_{-}{\bf 1}_{[0,t]}(y)
&=\frac{1}{\Gamma(H-\frac{1}{\alpha})}\int {\bf 1}_{[0,t]}(s)(s-y)_{+}^{H-\frac{1}{\alpha}-1}\e^{-\lambda(s-y)_{+}}\d s\\
&=\frac{1}{\Gamma(H-\frac{1}{\alpha}+1)}\Big[(t-y)_{+}^{H-\frac{1}{\alpha}}\e^{-\lambda(t-y)_{+}}-(-y)_{+}^{H-\frac{1}{\alpha}}\e^{-\lambda(-y)_{+}}
+ \lambda\int_{0}^{t}(s-y)_{+}^{H-\frac{1}{\alpha}}\e^{-\lambda(s-y)_{+}} \d s\Big]\\
&=\frac{1}{\Gamma(H-\frac{1}{\alpha}+1)}\ h_{H,\alpha,\lambda}(t;y),
\end{split}
\end{equation*}
see \eqref{hdef0}. This
proves \eqref{I1} and  \eqref{D1} follows similarly.
\hfill $\Box$

\begin{rem} \label{remlevy} {\rm
Note for $H = 1/\alpha$ \, $Z^{I\!I}_{1/\alpha,\alpha,\lambda} = M_\alpha$ is  $\alpha$-stable L\'evy process for  any $\lambda \ge 0$ which follows
from \eqref{eq:connection with fractional calculus}  and also  from \eqref{TFSM2}
by exchanging the order of integration:
\begin{eqnarray}\label{TFSM0}
Z^{I\!I}_{1/\alpha,\alpha,\lambda}(t)
&=&\int_{-\infty}^t \e^{-\lambda (t-y)}  M_\alpha (\d y)   - \int_{-\infty}^0 \e^{-\lambda (-y)_+} M_\alpha (\d y)  \\
&+&\lambda \int_{0}^{t} \d s  \int_{-\infty}^s  \e^{-\lambda(s-y)}\, M_\alpha (\d y)
= \int_0^t M_\alpha(\d y) =  M_\alpha(t). \nn
\end{eqnarray}
On the other hand,
\begin{equation}\label{TFSM00}
Z_{1/\alpha,\alpha,\lambda}(t) = Y_{\alpha,\lambda}(t) - Y_{\alpha,\lambda}(0)
\end{equation}
where $Y_{\alpha,\lambda}(t) := \int_{-\infty}^t \e^{-\lambda (t-y)}  M_\alpha (\d y), t \in \R$ is stationary
Ornstein-Uhlenbeck process.}
\end{rem}

\medskip

The next proposition gives the spectral domain representation of TFBM II.

\begin{prop}\label{prop:TFBMIIharmo}
Let $H>0, \lambda > 0$. The TFBM II \
$B^{I\!I}_{H,\lambda}$ in \eqref{TFBM2}
has the spectral domain representation
\begin{equation}\label{eq:har}
B^{I\!I}_{H,\lambda}(t)\eqfdd
\frac{1}{\sqrt{2\pi}}\int \frac{\e^{\i\omega t}-1}{\i\omega}(\lambda+\i\omega)^{\frac{1}{2}-H}\widehat{B}(\d\omega),
\end{equation}
where $\widehat{B}$ is an even complex-valued Gaussian white noise, $\overline {\widehat {B}(\d x)} = {\widehat B}(-\d x), $
with  zero  mean and variance $\E |{\widehat B}(\d x)|^2 = \d x$.
\end{prop}

\noi {\it Proof.}
From \eqref{I1}, $h_{H,2, \lambda}(t;y)=\Gamma (H + \frac{1}{2} ){\mathbb I}^{H - \frac{1}{2},\lambda}_{-} {\bf 1}_{[0,t]}(y)$ for $H>\frac{1}{2}$ and then
\begin{equation*}
\widehat{h_{H,2,\lambda}}(t;\omega) = \Gamma(H+\frac{1}{2})
\mathcal{F} \big[ {\mathbb I}^{H - \frac{1}{2},\lambda}_{-} {\bf 1}_{ [0,t]} \big](\omega)
=\frac{ \Gamma( H+\frac{1}{2} ) }{ \sqrt{2\pi}}\frac{\e^{\i\omega t}-1}{\i\omega}(\lambda+\i\omega)^{\frac{1}{2}-H},
\end{equation*}
where we used the Fourier transform of tempered fractional integrals (see Lemma 2.6 in \cite{MeerschaertsabzikarSPA}).
In the case $0<H<\frac{1}{2}$ the same expression for  $\widehat{h_{H,2,\lambda}}(t;\omega)$
follows by a similar argument.
Then \eqref{eq:har} follows from Parseval's formula for stochastic integrals, viz.,
$\int h_{H,2,\lambda}(t;y) B(\d y)\eqd \int \widehat{h_{H,2,\lambda}}(t;\omega){\widehat{B}}(\d\omega) $,
see Proposition 7.2.7
in \cite{SamorodnitskyTaqqu}.
\hfill $\Box$


\begin{rem}\label{rem:differentprocesses}
{\rm Meerschaert and Sabzikar \cite{Meerschaertsabzikar} obtained spectral representation of TFBM:
\begin{equation}\label{eq:TFBMdef}
B_{H,\lambda}(t)=\frac{\Gamma(H+\frac{1}{2})}{\sqrt{2\pi}}\int \frac{\e^{\i \omega t}-1}{(\lambda+
\i \omega)^{H+\frac{1}{2}}}\widehat{B}(\d\omega)
\end{equation}
where $H>0, \lambda > 0$.
By comparing spectral densities,
it can be shown that for $\lambda >0$ and $\sigma>0$, $B_{H,\lambda}$ and $\sigma B^{I\!I}_{H,\lambda}$ are different
processes (indeed, the coincidence of these spectral densities
would imply that $|(\lambda + \i \omega)/\omega |^2 $ is a constant function of $\omega $
which is possible if and only if $\lambda = 0$.) }

\end{rem}

The next proposition summarizes basic properties of TFSM II \, $Z^{I\!I}_{H,\alpha,\lambda}$.

\begin{prop}\label{lem:g squar integrable}

\smallskip

\noi (i) $Z^{I\!I}_{H,\alpha,\lambda}$
in \eqref{TFSM2} is well-defined for any $t \ge 0 $ and
$1< \alpha \le 2, H >0, \lambda > 0$,
as a stochastic integral in \eqref{Ialpha}.

\smallskip

\noi (ii) $Z^{I\!I}_{H,\alpha,\lambda}$
in \eqref{TFSM2} has stationary increments and
$\alpha$-stable finite-dimensional distributions. Moreover, it satisfies the following
scaling property:
\begin{equation}\label{eq:scalingTHP}
\left\{Z^{I\!I}_{H, \alpha, \lambda}(bt)\right\}_{t\in\R}\eqfdd
\left\{b^{H}Z^{I\!I}_{H,\alpha, b\lambda}(t)\right\}_{t\in\R}, \qquad \forall \ b >0.
\end{equation}

\smallskip

\noi (iii) $Z^{I\!I}_{H,\alpha,\lambda}$
in \eqref{TFSM2} has
a.s. continuous paths if either $\alpha = 2,  H>0$, or $ 1< \alpha < 2,
H > 1/\alpha$ hold.

\smallskip

\noi (iv)  The variance and covariance of TFBM II \,
$B^{I\!I}_{H,\lambda} $   $(H >0, \lambda >0)$ 
has the form
\begin{equation}\label{eq:variance}
\begin{split}
C^{2}_{t}&=\E \Big[(B^{I\!I}_{H,\lambda}(t))^2\Big]=\frac{1}{2\pi}\int \Big| \frac{\e^{\i\omega t}-1}{\i\omega}(\lambda+\i\omega)^{\frac{1}{2}-H}\Big|^{2}\ \d\omega\\
&=\frac{-2\Gamma(H)\lambda^{-2H}}{\sqrt{\pi}\Gamma(H- 1/2)}
\Big[1-{_2F_3}{ \Big(\{1,-1/2\}, \{1-H,1/2,1\},  \lambda^2 t^2/4\Big)}\Big]\\
&+t^{2H} \frac{\Gamma (1-H)}{ \sqrt{\pi} H 2^{2H} \Gamma (H+ 1/2)}\,{_2F_3} \Big(\{1,H- 1/2\},
\{1,H+1,H+ 1/2\}, \lambda^2 t^2/4\Big),
\end{split}
\end{equation}
and
\begin{equation}\label{eq:TFBM II acvf}
{\rm Cov}\Big[B^{I\!I}_{H,\lambda}(t),B^{I\!I}_{H,\lambda}(s)\Big]=\frac{1}{2}\Big[C_t^{2}+C_s^{2}-C_{t-s}^{2}\Big], \qquad
s,t\in\R,
\end{equation}
where $C_t^{2}$ is given in \eqref{eq:variance} and ${_2F_3}$ is the generalized hypergeometric function.
In particular,
\begin{eqnarray} \label{varBII}
\E (B^{I\!I}_{H,0}(t))^2
&=&\frac{t^{2H}}{\Gamma^2(H + 1/2)} \int_{-\infty}^t \big((t-s)^{H-1/2}- (-s)_+^{H-1/2}\big)^2 \d s \\
&=&t^{2H} \frac{\Gamma (1-H)}{\sqrt{\pi} H 2^{2H} \Gamma (H+ 1/2)}, \qquad 0< H < 1. \nn
\end{eqnarray}

\end{prop}

\noi {\it Proof.} (i) Follows from $h_{H, \alpha, \lambda}(t;\cdot) \in L^\alpha (\R)$, see above, also
 \cite{Meerschaertsabzikar2}.

\smallskip

\noi (ii) Stationarity of increments follows from the invariance properties
$h_{H, \alpha, \lambda}(t+T;y) - h_{H, \alpha, \lambda}(T;y)  = h_{H, \alpha, \lambda}(t; y-T) $ and
$\{ M_\alpha (y + T) - M_\alpha(T) \} \eqfdd \{M_\alpha (y)\},  \ \forall T \ \ge 0. $  Similarly,
property  \eqref{eq:scalingTHP}
follows from the scaling properties $h_{H, \alpha, \lambda/b}(bt;y) = b^{H- \frac{1}{\alpha}} h_{H, \alpha,\lambda}(t;y) $ and
$\{ M_\alpha (b t)\} \eqfdd \{ b^{\frac{1}{\alpha}} M(t) \}, \, \forall \, b>0. $

\smallskip

\noi (iii) We use the Kolmogorov criterion, see (\cite{Bill}, Theorem 12.4). Since $\lambda >0$ is fixed,
we can assume $\lambda =1 $  w.l.g.
First, let $\alpha = 2,  H >0$. Then $Z^{I\!I}_{H,\alpha,1} = C B^{I\!I}_{H,1}$ is a Gaussian process
with stationary increments.
Accordingly, it suffices to prove
$\E |B^{I\!I}_{H,1}(t)|^p  \le C t^\gamma $ for some $ p> 0, \gamma > 1$ and all $ 0 < t < 1 $. By Gaussianity,
$\E |B^{I\!I}_{H,1}(t)|^p \le C (\E |B^{I\!I}_{H,1}(t)|^2 )^{p/2} $.  We have
\begin{equation*}
\E |B^{I\!I}_{H,1}(t)|^2  = C \int_{-\infty}^t h^2_{H,2,1}(t;y) \d y =  C (I_1 + I_2),
\end{equation*}
where
\begin{equation*}
I_1 = \int_{-t}^t h^2_{H,2,1}(t;y) \d y \le C\int_0^{2t} (t - y)^{2(H- \frac{1}{2})} \d y
+ C t^2 (\int_0^{2t} s^{H - \frac{1}{2}} \d s)^2  \le C t^{2H}
\end{equation*}
and
\begin{equation*}
\begin{split}
I_2 = \int_{-\infty}^{-t} h^2_{H,2,1}(t;y) \d y &\le C\int_t^\infty ((t + y)^{H- \frac{1}{2}} \e^{-t-y}
- y^{H- \frac{1}{2}} \e^{-y})^2  \d y\\
&\  +  C\int_t^\infty (\int_0^t (s+y)^{H - \frac{1}{2}} \e^{-s-y} \d s)^2 \d y
=  C(I_2' + I_2'').
\end{split}
\end{equation*}
Using  $|(t + y)^{H- \frac{1}{2}} \e^{-t-y}
- y^{H- \frac{1}{2}} \e^{-y}| \le |\e^{-t}-1|\, \e^{-y} (t + y)^{H- \frac{1}{2}} +  \e^{-y}\, |(t + y)^{H- \frac{1}{2}}
- y^{H- \frac{1}{2}} | \le C t \, \e^{-y} (t + y)^{H- \frac{1}{2}} +  C t \, \e^{-y} y^{H- \frac{1}{2}} $ we obtain
$I'_2 \le C t^2 $ and, similarly, $I''_2 \le C t^2 $, implying
$I_1 + I_2 \le C(t^{2H} + t^2) $ and
$ \E |B^{I\!I}_{H,1}(t)|^p  \le C(t^{2H} + t^2)^{p/2}  $.
Hence, the above inequality is satisfied with $p > (1/H)\vee 1 $.
\\
Next, let $1 < \alpha < 2,  H > \frac{1}{\alpha} $.
Similarly as above, it suffices to prove
$\E |Z^{I\!I}_{H,\alpha,1}(t)|^p  \le C t^\gamma $ for some $ 1<  p < \alpha, \gamma > 1$ and all $ 0 < t < 1 $.
According to well-known moment inequality,
$\E |Z^{I\!I}_{H,\alpha,1}(t)|^p  = C \int_{-\infty}^t |h_{H,\alpha,1}(t;y)|^p \d y =  C (I_1 + I_2), $
where $I_1 = \int_{-t}^t |h_{H,\alpha,1}(t;y)|^p \d y \le C\int_0^{2t} (t - y)^{p(H- \frac{1}{\alpha})} \d  y
+ C t^p (\int_0^{2t} s^{H - \frac{1}{\alpha}} \d s)^p  \le C t^{1 + p(H - \frac{1}{\alpha})}  $
where $\gamma  =  1 + p(H - \frac{1}{\alpha}) > 1 $ due to $H > \frac{1}{\alpha} $. The bound
$I_2 = \int_{-\infty}^{-t} |h_{H,\alpha,1}(t;y)|^p \d y   \le C t^p $ can proved similarly as  in the case $\alpha =2$ above.
This proves part (iii).

\smallskip

\noi (iv) \eqref{eq:variance} follows from the spectral representation in \eqref{eq:har} and
the formula for integral transform in
(Prudnikov et al. \cite{Prudnikov}, p.379), see also
(\cite{Anh}, Lemma 2.1). In turn, \eqref{eq:TFBM II acvf} follows from  \eqref{eq:variance} and
stationarity of increments of $ B^{I\!I}_{H,\lambda}$.
Proposition  \ref{lem:g squar integrable} is proved.
\hfill $\Box$





\begin{rem} {\rm  (i)
Property \eqref{eq:scalingTHP}, as well as other properties in
Proposition  \ref{lem:g squar integrable} (ii), are also shared by the TFSM in \eqref{TFSM}, see
(\cite{Meerschaertsabzikar2},  Prop.2.3). It is related to the scaling property
for two-parameter processes introduced in  \cite{Jorgensen}.

\medskip

\noi (ii) 
For $H> 1/2 $, the covariance function of TFBM II $B^{I\!I}_{H,\lambda}$ admits the integral
representation
\begin{equation}\label{covZ}
\E B^{I\!I}_{H,\lambda}(t)B^{I\!I}_{H,\lambda}(s) =
C(H,\lambda)\int_{0}^{t}\int_{0}^{s}
|u-v|^{H-1}K_{H-1}(\lambda|u-v|)\d v\, \d u,
\end{equation}
where $C(H,\lambda)=\frac{2}{\sqrt{\pi}\Gamma(H-\frac{1}{2})(2\lambda)^{H-1}}$ and $K_{H-1} $ is the modified
Bessel function of second kind.
Formula \eqref{covZ}
follows from  \cite[p.344]{Gradshteyn} using the representation of $h_{H,2,\lambda} $ in \eqref{hdef1}.
For $H>1$ the integrand in \eqref{covZ}, viz.,
\begin{equation}
\frac{1}{\sqrt{\pi}\Gamma(H-\frac{1}{2})(2\lambda)^{H-1}}|u-v|^{H-1} K_{H-1}(\lambda|u-v|)
\end{equation}
is the  Mat\'{e}rn covariance function (in one dimension) with shape parameter $\nu=H-1>0$, scale parameter $\lambda>0$, and variance parameter 1, see e.g.  (\cite{Bolin}, (1.1)).
Note that the integral in \eqref{covZ} diverges when $0< H < 1/2$. A related
integral albeit more complex
representation of
the covariance
function of $B^{I\!I}_{H,\lambda}$ can be obtained for $0< H < 1/2 $, too,  
but we do not include it in the present paper. }
\end{rem}

The following theorem discusses local and global scaling properties
of TFSM and TFSM II.

\begin{thm}\label{lem:scaling}
{\it  Let $1< \alpha \le 2, 0< H<1 $ and $\lambda >0 $.

\smallskip

\noi (i) As $b \to \infty $
\begin{eqnarray}
b^{-1/\alpha}  Z^{I\!I}_{H, \alpha, \lambda}(bt) &\limfdd& c_{H,\alpha,\lambda}  M_\alpha(t) \quad \text{and} \quad
Z_{H, \alpha, \lambda}(bt) \ \limfdd\ C^+_{H,\alpha,\lambda} - C^-_{H,\alpha,\lambda},  \label{large}
\end{eqnarray}
where $c_{H,\alpha,\lambda} = \lambda^{\frac{1}{\alpha} - H} \Gamma(1 + H - \frac{1}{\alpha})$ and
$C^+_{H,\alpha,\lambda}, C^-_{H,\alpha,\lambda}$ are independent copies
of $\alpha$-stable random variable $C^0_{H,\alpha,\lambda} $ in \eqref{C0}.
\smallskip

\noi (ii) As  $b \to 0 $
\begin{eqnarray}
b^{-H}  Z^{I\!I}_{H, \alpha, \lambda}(bt) &\limfdd& c_{H,\alpha}  Z^{I\!I}_{H,\alpha,0}(t) \quad \text{and} \quad
b^{-H}  Z_{H, \alpha, \lambda}(bt) \ \limfdd\ c_{H,\alpha}  Z_{H,\alpha,0}(t), \label{small}
\end{eqnarray}
where $c_{H,\alpha} $ is defined in \eqref{clim2} below and $ Z^{I\!I}_{H,\alpha,0}  = Z_{H,\alpha,0}$
is fractional stable motion.
}
\end{thm}

\noi {\it Proof.} For simplicity of notation we shall assume that $\beta = 0$ ($M_\alpha$ is symmetric)  and
$\sigma = 1$ in \eqref{Mcf}. The proof in the general case is analogous.

\smallskip

\noi (i) Consider the first relation in \eqref{large}.
It suffices to prove the convergence of
characteristic functions
$$
\E \exp\{ \i  b^{- \frac{1}{\alpha}}\sum_{j=1}^m \theta_j  (Z^{I\!I}_{H, \alpha, \lambda}(b t_j)  -
Z^{I\!I}_{H, \alpha, \lambda}(bt_{j-1}))   \}
\ \to \ \E \exp\{ \i \sum_{j=1}^m \theta_j  (M_{\alpha}(t_j) - M_\alpha(t_{j-1}))\}, \qquad b \to \infty
$$
for any $0= t_0 < t_1 < \dots < t_m, \theta_j \in \R, m=1,2,\dots $,  or
\begin{eqnarray} \label{clim1}
b^{-1} \int_{\R} \big|\sum_{j=1}^m \theta_j (h_{H,\alpha,\lambda} (bt_j;y) - h_{H,\alpha,\lambda} (bt_{j-1}; y)) \big|^\alpha \d y
&\to&\sum_{j=1}^m |\theta_j|^\alpha (t_j - t_{j-1})
\end{eqnarray}
as $b  \to \infty$.
Consider \eqref{clim1} for $m=1, t_1 = 1 $ (the proof in the  general case seems is similar). This follows from
\begin{equation}\label{hconv}
b^{-1} \int_0^{c} |h_{H, \alpha, \lambda}(b;y)|^\alpha \d y \to  c^\alpha_{H, \alpha,\lambda} \quad  \text{and} \quad
b^{-1} \int_{-\infty}^{0} |h_{H, \alpha,\lambda}(b; y)|^\alpha \d  y \to 0.
\end{equation}
Note for each $ y>0$
\begin{eqnarray}\label{hlim}
c_{H, \alpha,\lambda}&=&\lim_{b \to \infty} h_{H,\alpha,\lambda}(b;y)
= \lim_{b \to \infty}   \lambda \int_{0}^{b} (s-y)_{+}^{H-\frac{1}{\alpha}}\e^{-\lambda(s-y)_{+}}\ \d s.
\end{eqnarray}
Since $b^{-1} \int_0^{b} |h_{H, \alpha, \lambda}(c;y)|^\alpha \d y \le 2(1+ \lambda) \int_0^\infty z^{\alpha H -1} \e^{-\alpha \lambda z} \d z
 < \infty $
for all $b \ge 1$, the first  convergence in \eqref{hconv}  follows from \eqref{hlim} and the dominated convergence
theorem.  The second convergence  in \eqref{hconv}  follows by noting that
$\limsup_{b \to \infty } \int_{-\infty}^{0} |h_{H, \alpha,\lambda}(b; y)|^\alpha \d y  < \infty $ since
\begin{equation}
\begin{split}
\int_{-\infty}^{0} |h_{H, \alpha,\lambda}(b; y)|^\alpha \d y \ &\le \ 6 \int_0^\infty z^{\alpha H -1} \e^{-\lambda \alpha z} \d z
\ + \ 3 \lambda^\alpha
\int_0^\infty  \d y \big( \int_0^\infty (s+y)^{H- \frac{1}{\alpha}} \e^{-\lambda (s+y)} \d s \big)^\alpha\\
&=  6 \Gamma (\alpha H)/(\lambda \alpha)^{\alpha H} + 3 \lambda^\alpha
\int_0^\infty \d  y \big(\int_y^\infty w^{H- \frac{1}{\alpha}} \e^{-\lambda w} \d w \big)^\alpha
    < \infty
\end{split}
\end{equation}
This proves the first convergence in \eqref{large}.

Consider the second convergence in \eqref{large} for $t=1 $. It suffices to show that
\begin{equation}\label{gconv}
\int_{-\infty}^{b} \Big|g_{H, \alpha, \lambda}(b;y)\Big|^\alpha \d y \to 2\int_0^\infty y^{\alpha H -1}  \e^{-\alpha \lambda y} \d y
= \frac{2\Gamma(H\alpha)}{(\lambda\alpha)^{H\alpha}}
\end{equation}
as $b\to\infty$. Here,
\begin{equation}\label{eq:g kernel}
g_{H, \alpha, \lambda}(t;y) :=(t-y)_+^{H - \frac{1}{\alpha}} \e^{-\lambda (t-y)_+} - (-y)_+^{H - \frac{1}{\alpha}} \e^{-\lambda (-y)_+},
\qquad y \in \R
\end{equation}
is the kernel of TFSM \eqref{TFSM}. Relation  \eqref{gconv}
follows from $\int_0^{b} \Big|g_{H, \alpha, \lambda}(b;y)\Big|^\alpha \d y \to \Gamma(H\alpha)/(\lambda\alpha)^{H\alpha}$ and
$\int_{-\infty}^0 \Big|g_{H, \alpha, \lambda}(b;y)\Big|^\alpha \d y \to  \Gamma(H\alpha)/(\lambda\alpha)^{H\alpha}$
using the dominated convergence theorem.
The general finite-dimensional convergence
follows analogously.
This  proves part (i).

\smallskip

\noi (ii)
As in part (i), we restrict the proof of  \eqref{small}
to one-dimensional convergence at $t=1$ since the general case
follows analogously.  Consider the first convergence in  \eqref{small}. It suffices to show
\begin{eqnarray} \label{clim2}
b^{-\alpha H} \int_{\R} \big|h_{H, \alpha,\lambda} (b;y)\big|^\alpha \d y
&\to&c^\alpha_{H,\alpha} = \int_{\R} |h_{H,\alpha, 0} (1;y)|^\alpha \d y
\end{eqnarray}
as $b \to 0$. By change of variables, 
we have $b^{-\alpha H} \int_{\R} \big|h_{H, \alpha,\lambda} (b;y)\big|^\alpha \d y
= \int_{\R} \big|h_{H, \alpha, b\lambda} (1;y)\big|^\alpha \d y $ where
$h_{H, \alpha, b\lambda} (1;y) \to  h_{H, \alpha, 0} (1;y) \ (b \to 0) $ for
each $y \in \R,  y \ne 1 $ and the convergence
of integrals in \eqref{clim2} can be justified by the dominated convergence theorem.
The second convergence in  \eqref{small} for $t=1$ follows similarly from
$b^{-\alpha H} \int_{\R} \big|g_{H, \alpha,\lambda} (b;y)\big|^\alpha \d y
\to c^\alpha_{H,\alpha}$ and we omit the details.
Theorem  \ref{lem:scaling} is proved.  \hfill $\Box $

\vskip.3cm

\section{Dependence properties of TFSN II}

Recall the definition of tempered fractional stable noise
(TFSN II)  $Y^{I\!I}_{H, \alpha, \lambda}$
in \eqref{eq:TFSN definition}, which is a stationary process with discrete time $t\in \Z$.
The following proposition obtains the spectral representation and spectral density of this process, or
TFGN II
\begin{equation}\label{eq:TFGNIIdef}
Y^{I\!I}_{H,\lambda} =
\big\{Y^{I\!I}_{H,\lambda}(t) :=  B^{I\!I}_{H,\lambda}(t+1) - B^{I\!I}_{H,\lambda}(t), \ t \in \Z\big\}
\end{equation}
in the Gaussian case $\alpha = 2$.

\begin{prop}\label{spectraldesnityTFGNII}
TFGN II in \eqref{eq:TFGNIIharmonizable} has spectral representation
\begin{equation}\label{eq:TFGNIIharmonizable}
Y^{I\!I}_{H,\lambda}(t) =
\frac{1}{\sqrt{2\pi}} \int \e^{\i \omega t} \frac{\e^{\i\omega}-1}{\i\omega}({\lambda+ \i\omega})^{\frac{1}{2}-H}\widehat{B}(\d \omega),
\quad t \in \Z
\end{equation}
and spectral density
\begin{equation}\label{eq:spectraldesnityTFGNII}
h(\omega)= \frac{1}{2\pi} \Big\{ \big|\frac{\e^{\i\omega}-1}{\i\omega}\big|^2 [\lambda^2+\omega^2]^{\frac{1}{2}-H} +
|\e^{\i\omega}-1|^2
\sum_{\ell \in \Z \setminus \{0\}} \frac{1}{(\omega + 2\pi \ell)^2
[\lambda^2+(\omega+2\pi \ell)^2]^{H-\frac{1}{2}}}\Big\}, \qquad \omega \in [-\pi, \pi].
\end{equation}
\end{prop}
{\it Proof.} \eqref{eq:TFGNIIharmonizable} is immediate from \eqref{eq:har}. Whence it follows that the covariance
\begin{equation}\label{specdenscalc}
\begin{split}
r(j)
=&\frac{1}{2\pi}\int \e^{\i\omega j}\Big|\frac{\e^{\i\omega}-1}{\i\omega}\Big|^{2}(\lambda^2+\omega^2)^{\frac{1}{2}-H}\d\omega\\
=&\frac{1}{2\pi}\int^{+\pi}_{-\pi}\e^{\i\omega j}\big|\e^{\i\omega}-1\big|^2
\sum_{\ell\in \Z} (\omega +2\pi \ell)^{-2}  {[\lambda^2+(\omega+2\pi \ell)^2]}^{\frac{1}{2}-H}\d\omega, \qquad j \in \Z
\end{split}
\end{equation}
implying  \eqref{eq:spectraldesnityTFGNII} and the proposition. \hfill $\Box$

\smallskip

Note that the spectral density in \eqref{eq:spectraldesnityTFGNII} is bounded, continuous and separated from zero on the whole interval
$[-\pi, \pi]$. Particularly, for any $H, \lambda>0$
\begin{equation*}
h(0) = \lim_{\omega \to 0} h(\omega) =   \lambda^{1-2H}/2\pi  \in (0, \infty).
\end{equation*}
Using the popular terminology (see e.g. \cite{koul}, Definition 3.1.4) we may say that
TFGN II $Y^{I\!I}_{H,\lambda}$ has short memory in the frequency domain for any $H, \lambda >0$ while
TFGN has negative memory in the frequency domain since its spectral density vanishes at $\omega =0$.
The distinctions
in low frequency spectrum  between the two processes  are illustrated in
Figures 1 and 2. 

\begin{figure}[ht]\label{fig:01}
\includegraphics[scale=0.5]{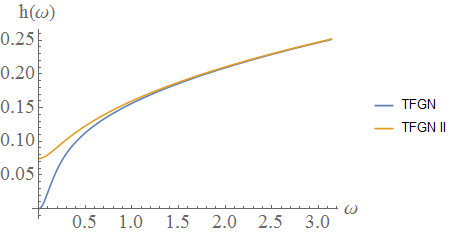}
\includegraphics[scale=0.5]{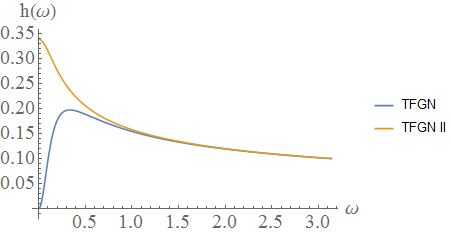}
\caption{Left: spectral density of TFGN and TFGN II with $H=0.3$ and $\lambda=0.15$. Right: the same plot for $H=0.7$ and $\lambda=0.15$.}
\end{figure}

\begin{figure}[ht]\label{fig:03}
\includegraphics[scale=0.5]{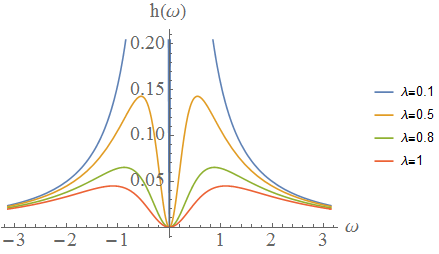}
\includegraphics[scale=0.5]{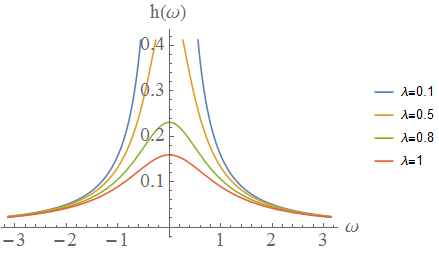}
\caption{Left: spectral density of  TFGN with $H=4/3$ and $\lambda=0.1, 0.5, 0.8,$ and $1$. Right:
spectral density of TFGN II for the same values of H and $\lambda$}.
\end{figure}


\begin{rem} \label{kolm}
{\rm In the celebrated Kolmogorov's  model  of turbulence  (\cite{Kolmogorov}, \cite{Kolmogorov1}),
the spectral density of turbulent velocity data in the
inertial range is proportional to $\omega^{-5/3}$ for moderate frequencies $\omega$ and is
bounded at low frequencies. Such behavior can be exhibited by both TFGN and
TFGN II with $H =4/3$, see  \eqref{eq:spectraldesnityTFGNII}.
Meerschaert et al.  \cite{Meerschaertsabzikarkumarzeleki}  validate
ARTFIMA$(0,d,\lambda,0)$ model on turbulent water velocities in the Great Lakes region.
They find that this model effectively captures both the correlation properties
and the underlying probability distribution of this data, with $d = 5/6$ corresponding to Kolmogorov's
law, and small tempering parameter $\lambda = 0.006$. Since ARTFIMA$(0,d,\lambda,0)$ is closely
related to TFGN II, see \cite{SabzikarSurgailis}, we expect that
the latter model can be also successfully applied for modeling of turbulent data.

}
\end{rem}


Next, we discuss dependence properties of TFSN II with $ 1< \alpha < 2 $.
Since  this process 
has infinite variance, 
other numerical characteristics extending the notion of covariance
must be used to characterize the decay rate of dependence between  $Y^{I\!I}_{H,\alpha,\lambda}(t_1)$ and
$Y^{I\!I}_{H,\alpha,\lambda}(t_2)$ as $t_2 - t_1 \to \infty $.  Astrauskas et al.\ \cite{Astrauskas},
Samorodnitsky and Taqqu \cite{SamorodnitskyTaqqu}, Koul and Surgailis \cite{koulsurgailis}
discussed the asymptotics as $t \to \infty $ of the bivariate characteristic function
\begin{equation}\label{eq: r definition}
r_Y(t;\theta_1,\theta_2) = \mathbb{E}\big[\e^{\i(\theta_1 Y(t) +\theta_2 Y(0))}\big] -\mathbb{E}
\big[\e^{\i\theta_1 Y(t)}\big] \mathbb{E}\big[\e^{\i\theta_2Y(0)}\big], \quad \theta_1, \theta_2 \in\mathbb{R}
\end{equation}
for some stationary  processes $Y = \{Y(t), t \in \Z\} $
with infinite variance, including fractional stable noise (FGN) and long memory moving average with infinite variance
innovations. Meerschaert and Sabzikar \cite{Meerschaertsabzikar2} studied the asymptotic behavior of \eqref{eq: r definition} for
TFSN $Y_{H,\alpha,\lambda}(t) =  Z_{H,\alpha,\lambda}(t+1)-Z_{H,\alpha,\lambda}(t),  t \in \Z $,
where $Z_{H,\alpha,\lambda}$ is the TFSM in \eqref{TFSM}.
Given two real-valued functions $f(t)$, $g(t)$ on $\mathbb{R}$, we will write $f(t)\asymp g(t)$ if $C_1\leq |{f(t)}/{g(t)}|\leq C_2$ for all $t>0$ sufficiently large, for some $0<C_1<C_2<\infty$.
(\cite{Meerschaertsabzikar2}, Theorem 2.6) proved that for $1<\alpha<2$, $\frac{1}{\alpha}<H<1$
\begin{equation}\label{asymptotic r for TFSN I}
r_{Y_{H,\alpha,\lambda}}(t; \theta_1, \theta_2) \asymp \e^{-\lambda t} t^{H-\frac{1}{\alpha}}, \quad\ t\to\infty.
\end{equation}
Below, we prove the following result about the behavior  of \eqref{eq: r definition} for TFSN II which is another indication 
that TFSM and TFSM II are  different processes.

\begin{thm}\label{thm:asym1} Let  $Z^{I\!I}_{H,\alpha,\lambda}$ be a TFSM II in \eqref{TFSM2} with
S$\alpha$S noise $M_\alpha $,  $1<\alpha \le 2, \, \frac{1}{\alpha} < H <1,  \,   \lambda >0$, and
$Y^{I\!I}_{H,\alpha,\lambda}(t)$ be TFSN II in  \eqref{eq:TFSN definition}.  Then for any
$\theta_1, \theta_2 \ne 0 $
\begin{equation}\label{mmm1}
r_{Y^{I\!I}_{H,\alpha,\lambda}}(t;\theta_1, \theta_2) \asymp    \e^{-\lambda t}   
t^{H-\frac{1}{\alpha}-1}, \quad t \to \infty.
\end{equation}
\end{thm}

\noi {\it Proof.} For concreteness, assume $\theta_1 \theta_2 >0 $, the case $\theta_1 \theta_2 <0$ being analogous.
Write $r(t) =  r_{Y^{I\!I}_{H,\alpha,\lambda}}(t;\theta_1, \theta_2), Y(t) = Y^{I\!I}_{H,\alpha,\lambda}(t),
h(t;x) = h_{H,\alpha, \lambda}(t;x) $ for brevity.
Then
$r(t) =  K \big(\e^{-I(t)}-1\big) $,
where $K = K(\theta_1,\theta_2) =
\E [\e^{\i\theta_1Y(0)}]
\E [\e^{\i\theta_2Y(0)}] >0 $ and
\begin{equation}\label{It}
I(t) = I(t;\theta_1,\theta_2):= \|\theta_1Y(t)+\theta_2Y(0)\|_\alpha^\alpha-\|\theta_1Y(t)\|_\alpha^\alpha-\|\theta_2Y(0)\|_\alpha^\alpha.
\end{equation}
Clearly, it suffices to prove \eqref{mmm1} for $I(t)$, viz.,
\begin{equation}\label{mmm2}
I(t) \asymp    \e^{-\lambda t} t^{H-\frac{1}{\alpha}-1}, \quad t \to \infty.
\end{equation}
It follows from \eqref{TFSM2} that  $Y(t) = \int_{\R} g(t;x) M_\alpha (\d x) $ where
$g(t;x) = h(t+1;x)- h(t;x)$. From the representation  \eqref{hdef1} we obtain
\begin{equation}\label{gdef}
g(t;x) =  (H- \frac{1}{\alpha}) \int_t^{t+1} (s-x)_+^{H-\frac{1}{\alpha}-1} \e^{-\lambda(s-x)_{+}} \d  s, \   -\infty < x < t.
\end{equation}
We have
$$
I(t)  = \int_{-\infty}^1 \big\{ |\theta_1 g(t;x) + \theta_2 g(0;x)|^\alpha
- |\theta_1 g(t;x)|^\alpha - |\theta_2 g(0;x)|^\alpha \big\}\d x
$$
Then $g(t;x) >0, \, x < t+ 1$ and for any $-\infty < x < 1  < t$ we obtain
\begin{eqnarray*}
\frac{g(t;x)}{g(0;x)} = \frac{\int_t^{t+1} (s-x)^{H-\frac{1}{\alpha}-1} \e^{-\lambda(s-x)} \d s}
{\int_0^{1} (s-x)_+^{H-\frac{1}{\alpha}-1} \e^{-\lambda(s-x)} \d s}
= \frac{\int_t^{t+1} (s-x)^{H-\frac{1}{\alpha}-1} \e^{-\lambda s} \d s}
{\int_0^{1} (s-x)^{H-\frac{1}{\alpha}-1} \e^{-\lambda s} \d s} \le \e^{-\lambda (t-1)} \frac{g^0(t;x)}{g^0(0;x)}
\end{eqnarray*}
where
$$
g^0(t;x) :=   (H- \frac{1}{\alpha}) \int_t^{t+1} (s-x)_+^{H-\frac{1}{\alpha}-1} \d s =  (t+1 -x)_+^{H-\frac{1}{\alpha}} - (t -x)_+^{H-\frac{1}{\alpha}}.  $$
Similarly,
\begin{eqnarray*}
\frac{g(t;x)}{g(0;x)}
= \frac{\int_t^{t+1} (s-x)^{H-\frac{1}{\alpha}-1} \e^{-\lambda s} \d s}
{\int_0^{1} (s-x)^{H-\frac{1}{\alpha}-1} \e^{-\lambda s}\d s} \ge \e^{-\lambda (t+1)} \frac{g^0(t;x)}{g^0(0;x)}
\end{eqnarray*}
implying
\begin{equation}\label{gdomin}
C_1 \e^{-\lambda  t} \frac{g^0(t;x)}{g^0(0;x)} \le   \frac{g(t;x)}{g(0;x)} \le  C_2  \e^{-\lambda t} \frac{g^0(t;x)}{g^0(0;x)}
\end{equation}
for some constants $0 < C_1 < C_2 < \infty $ independent of $-\infty < x < 1  < t$.
Then $I(t) = I_1(t) - I_2(t)$, where
\begin{eqnarray*}
I_1(t)&:=& \int_{-\infty}^1  \big|\theta_2 g(0;x)\big|^\alpha \Big\{ \big| 1+  \frac{\theta_1}{\theta_2} \frac{ g(t;x)}{g(0;x)}\big|^\alpha
- 1 \Big\} \d x, \qquad
I_2(t) :=  \int_{-\infty}^1 \big|\theta_1g(t;x)\big|^\alpha  \d x.
\end{eqnarray*}
Next, split $I_1(t) = \int_{-\infty}^0 \cdots + \int_{0}^1 \cdots =: I_1'(t) + I_1''(t)$. Note $g^0(t;x) \le g^0(0;x)$ for
$x <0, t >0  $ and therefore  $\sup_{x <0} \frac{g(t;x)}{g(0;x)} \le C_2 \e^{-\lambda t} \to 0 $ as $t \to \infty $.
Therefore and from  \eqref{gdomin}  we obtain
for $x <0$
\begin{equation}
C_1 \e^{-\lambda  t} \frac{g^0(t;x)}{g^0(0;x)} \le
C_1  \frac{ g(t;x)}{g(0;x)} \le   \big| 1+  \frac{\theta_1}{\theta_2} \frac{ g(t;x)}{g(0;x)}\big|^\alpha - 1  \le
C_2 \frac{ g(t;x)}{g(0;x)} \le  C_2 \e^{-\lambda  t} \frac{g^0(t;x)}{g^0(0;x)}.
\end{equation}
This implies
\begin{equation}\label{I1bdd}
I'_1(t) \asymp  \e^{-\lambda t} J(t), \quad \text{where} \quad J(t) = \int_{-\infty}^0  \frac{ g(0;x)^\alpha }{g^0(0;x)} \, g^0(t,x)\d  x.
\end{equation}
Let us prove that
\begin{equation}\label{Jbdd}
J(t) \asymp t^{H-\frac{1}{\alpha}-1}.
\end{equation}
Split $J(t)  =  \int_{-\infty}^{-2} \cdots + \int_{-2}^0 \cdots =:
J_1(t) + J_2(t)$. Then since $g(0;x) \le g^0(0;x), x <0$ we have that
$J_2(t) \le C \int_0^2 ((t+1 +x)^{H- \frac{1}{\alpha}} - (t+x)^{H - \frac{1}{\alpha}})\d x  \le  C t^{H - \frac{1}{\alpha}-1} $.
Next, using $C_1z^{H -\frac{1}{\alpha}-1}  \le  (z+1)^{H - \frac{1}{\alpha}} - z^{H - \frac{1}{\alpha}} \le
C_2 z^{H-   \frac{1}{\alpha} -1}, z \ge 1 $
we obtain
\begin{eqnarray*}
J_1(t)&\asymp&\int_2^\infty \e^{-\alpha \lambda x} \big( (1+x)^{H - \frac{1}{\alpha}} - x^{H - \frac{1}{\alpha}}\big)^{\alpha -1}
 \big( (t+1+x)^{H - \frac{1}{\alpha}} - (t+x)^{H - \frac{1}{\alpha}}\big)\d  x \\
&\asymp&\int_2^\infty \e^{-\alpha \lambda x} x^{H -\frac{1}{\alpha}-1}  (t+x)^{H -\frac{1}{\alpha}-1} \d  x \\
&\asymp&t^{H -\frac{1}{\alpha}-1}
\end{eqnarray*}
proving \eqref{Jbdd}. Next, let us show that
\begin{equation}\label{I2bdd}
I''_1(t) \le C \e^{-\lambda t}  t^{H-\frac{1}{\alpha}-1}
\end{equation}
Indeed, from  \eqref{gdomin} and the definition of $g(0;x) $ and $g^0(t;x)$ it follows that
$$
I''_1(t) \le C  \int_0^1 (1-x)^{(H-\frac{1}{\alpha}-1)\alpha} \Big\{ \big(1 +  C \frac{\e^{-\lambda t} t^{H-\frac{1}{\alpha}-1}}{(1-x)^{H-\frac{1}{\alpha}}} \big)^\alpha - 1 \Big\}\d x \le
C \e^{-\lambda t}  t^{H-\frac{1}{\alpha}-1},
$$
proving  \eqref{I2bdd}. Relation
\begin{equation}\label{I3bdd}
I_1(t) \asymp \e^{-\lambda t}  t^{H-\frac{1}{\alpha}-1}
\end{equation}
follows from \eqref{I1bdd}, \eqref{Jbdd}, \eqref{I2bdd} and
the fact that $I'_1(t) \ge 0, I''_2(t) \ge 0$, due to $\theta_1 \theta_2 >0$. \\
It remains to show
\begin{equation} \label{I4bdd}
I_2(t) =  o\big(\e^{-\lambda t}  t^{H-\frac{1}{\alpha}-1}\big).
\end{equation}
Split $I_2(t) =  \int_{-\infty}^0 \cdots + \int_0^1 \cdots  = I'_2(t) + I''_2(t)$.
According to \eqref{gdomin}, $I'_2(t) \le C \e^{-\alpha \lambda t} \tilde J(t) $ where
$$
\tilde J(t) = \int_{-\infty}^0    \big( \frac{ g^0(t;x) }{g^0(0;x)} \big)^\alpha \, g(0;x)^\alpha \d  x
\le J(t)
$$
since $g^0(t;x) \le  g^0(0;x),  x < 0$, with $J(t)$ as in
\eqref{I1bdd}, \eqref{Jbdd}. Then, $I'_2(t)$ satisfies
\eqref{I4bdd} since $\alpha > 1 $. Finally, since $g(t;x) \le  C \e^{-\lambda t} t^{H-\frac{1}{\alpha}-1}, 0< x <1 $ so
$I''_2(t) \le  C \e^{-\alpha \lambda t}  t^{(H-\frac{1}{\alpha}-1)\alpha} =  o\big(\e^{-\lambda t}  t^{H-\frac{1}{\alpha}-1}\big)$,
ending the proof of \eqref{I4bdd} and \eqref{mmm2}.  Theorem
\ref{thm:asym1} is proved.  \hfill  $\Box$

\begin{rem} 
{\rm  Note that in the case $\lambda=0$,  the exponent $H - \frac{1}{\alpha} - 1$ in
\eqref{mmm1} does not agree with the corresponding exponent for (untempered) fractional stable noise
in Astrauskas et al.\ \cite{Astrauskas}, including the Gaussian case $\alpha =2$. Particularly, for fractional Gaussian noise
the covariance function as well as the
function $r_Y(t)$ in \eqref{eq: r definition} decay at rate $t^{2H-2} $, see  \cite{SamorodnitskyTaqqu}.
The reason is that \eqref{mmm1} holds for {\it fixed} $\lambda >0$, or  the tempered fractional noise
alone. 
}
\end{rem}

\begin{rem} 
{\rm  In the Gaussian case $\alpha = 2$ we have that \eqref{It} is proportional to the covariance
of TFGN II:  $I(1,1;t) = {\rm Cov}(Y^{I\!I}_{H,\lambda}(0), Y^{I\!I}_{H,\lambda}(t)).$  Therefore
\eqref{mmm2}  implies
\begin{equation*}
{\rm Cov}(Y^{I\!I}_{H,\lambda}(0), Y^{I\!I}_{H,\lambda}(t)) \asymp    \e^{-\lambda t} t^{H-\frac{3}{2}}, \quad t \to \infty,
\end{equation*}
for  $\lambda >0,  1/2 < H <1.   $
}
\end{rem}


\bigskip


\footnotesize

\begin{thebibliography}{99}



\bibitem{Anh}
Anh, V., Heyde, C. and Tieng, Q. (1999).
Stochastic models for fractal processes. {\it J. Statist. Plann. Inference} {\bf 80}\ 123--135.




\bibitem{Astrauskas} Astrauskas, A., Levy, J.B.  and Taqqu, M.S. (1991).
The asymptotic dependence structure of the linear fractional L\'evy motion. {\it Lithuanian Math. J.} {\bf 31}\  1--19.






\bibitem{Bill}
Billingsley, P. (1968). {\it Convergence of Probability Measures}. New York: Wiley.


\bibitem{Bolin}
Bolin, D. and Lindgren, F. (2011).
Spatial models generated by nested stochastic partial differential equations, with an application to global ozone mapping.
{\it Ann. Appl. Statist.} {\bf 5}\   523--550.






\bibitem{koul}
Giraitis, L., Koul. H.L. and Surgailis, D. (2012). {\it Large Sample Inference for Long Memory Processes.}
London: Imperial College Press.


\bibitem{Gradshteyn}
Gradshteyn, I.S. and Ryzhik, I.M. (2000).
{\it Tables of Integrals and Products.} 6th edition. New York: Academic Press.


\bibitem{Jorgensen} J{\o}rgensen, B., Martinez, J.R. and Demetrio, C.G.B. (2011).
Self-similarity
and Lamperti convergence for families of stochastic processes.
{\it Lithuanian Math. J.} {\bf 51}\  342--361.


\bibitem{Kolmogorov}
Kolmogorov, A.N. (1940).
Wiener spiral and some other interesting curves in Hilbert space.  {\it Dokl. Akad. Nauk SSSR} {\bf 26}\ 115--118.

\bibitem{Kolmogorov1}
Kolmogorov, A.N. (1941).
The local structure of turbulence in an incompressible fluid at very high Reynolds numbers.
{\it Dokl. Akad. Nauk SSSR} {\bf 30}\ 209--303.


\bibitem{koulsurgailis}
Koul, H.L. and Surgailis, D. (2001). Asymptotics of empirical processes of long memory
moving averages with infinite variance. {\it  Stoch. Process.  Appl.}  {\bf 91} \ 309--336.



\bibitem{MeerschaertsabzikarSPA}
Meerschaert, M.M. and Sabzikar, F. (2014).
Stochastic integration with respect to tempered fractional Brownian motion.
{\it Stochastic Process. Appl.} {\bf 124} \ 2363--2387.



\bibitem{Meerschaertsabzikar}
Meerschaert, M.M. and Sabzikar, F. (2013).
Tempered fractional Brownian motion. {\it  Statist. Probab. Lett.} {\bf 83} \ 2269--2275.



\bibitem{Meerschaertsabzikar2} Meerschaert, M.M. and Sabzikar, F. (2016). Tempered fractional stable motion.
{\it J. Theoret. Probab.} {\bf 29}\ 681--706.


\bibitem{Meerschaertsabzikarkumarzeleki}
Meerschaert, M.M., Sabzikar, F., Phanikumar, M.S. and Zeleke, A. (2014).
Tempered fractional time series model for turbulence in geophysical flows.
{\it J. Stat. Mech. Theory  Exp.} {\bf 2014} \  P09023.







\bibitem{Prudnikov}
Prudnikov, A.P., Brychkov, Y.A. and Marichev, O.I. (1986). {\it  Integrals and Series, vol. 1.} New York: Gordon and Breach.

\bibitem{SamorodnitskyTaqqu}
Samorodnitsky, G. and Taqqu, M.S.  (1996). {\it Stable Non-Gaussian Random Processes: Stochastic Models with Infinite Variance.}
Boca Raton etc: Chapman and Hall.

\bibitem{SabzikarSurgailis} Sabzikar, F. and Surgailis, D. (2017).  Invariance  principles for tempered fractionally
integrated processes. Preprint.





















\end{thebibliography}
\end{document}